# Solution of Conformable Fractional Ordinary Differential Equations via Differential Transform Method

Emrah Ünal[a], Ahmet Gökdoğan[b]

[a] Department of Elementary Mathematics Education, Artvin Çoruh University, 08100 Artvin, Turkey

emrah.unal@artvin.edu.tr

[b] Department of Mathematical Engineering, Gümüşhane University, 29100 Gümüşhane, Turkey,

gokdogan@gumushane.edu.tr

**Abstract**

Recently, a new fractional derivative called the conformable fractional derivative is given which is based on the basic limit definition of the derivative in [1]. Then, the fractional versions of chain rules, exponential functions, Gronwall's inequality, integration by parts, Taylor power series expansions is developed in [2]. In this paper, we give conformable fractional differential transform method and its application to conformable fractional differential equations.

**Key words :** Conformable Fractional Derivative, Fractional power series, Conformable Fractional Differential Transform Method, Conformable Fractional Ordinary Differential Equations

## 1. Introduction

Despite becoming a popular topic in recent years, the concept of fractional derivatives has emerged in the late 17th century. Several definitions have been made to define the fractional derivative and continues to be done. Most popular definitions in this area are the Riemann-Liouville, Caputo and Grunwald-Letnikov definitions. These definitions are defined as, respectively,

(I)  Riemann Liouville definition:
$$D_x^\alpha f(x) = \frac{1}{\Gamma(n-\alpha)} \left(\frac{d}{dx}\right)^n \int_0^x (x-t)^{n-\alpha-1} f(t) dt, \qquad n-1 < \alpha \leq n$$

(II) Caputo definition:
$$D_x^\alpha f(x) = \frac{1}{\Gamma(n-\alpha)} \int_0^x (x-t)^{n-\alpha-1} f^{(n)}(t) dt, \qquad n-1 < \alpha \leq n$$

(III) Grunwald-Letnikov definition:
$$_aD_x^\alpha f(x) = \lim_{h\to 0} h^{-\alpha} \sum_{j\to 0}^{\frac{x-a}{h}} (-1)^j \binom{\alpha}{j} f(x-jh)$$

Recently, a new definition of fractional derivative [1] that prominently compatible with the classical derivative was made by Khalil et al. Unlike other definitions, this new definition satisfies formulas of derivative of product and quotient of two functions and has a simpler the chain rule. In addition to conformable fractional derivative definition, the conformable fractional integral definition, Rolle theorem and Mean value theorem for conformable fractional differentiable functions was given. Another study [2] was done by Abdeljawad that contributed to this new field. He presented left and right conformable fractional derivative, fractional integrals of higher orders concepts. Moreover, he gave the fractional chain rule, the fractional integration by parts formulas, Gronwall inequality, the fractional power series expansion and the fractional Laplace transform definition. In short time, many studies related to this new fractional derivative definition was done [3,4,5]

The differential transform method (DTM) is one of the numerical methods that is used for finding the solution of differential equations. The DTM was first proposed by Zhou. He solved linear and nonlinear initial value problems in electric circuit analysis in [6] via DTM. DTM is used in many studies related to the eigenvalue problems [7,8], the linear and non-linear higher-order boundary value problems [9], the higher-order initial value problems [10,11], systems of ordinary and partial differential equations [12,13,14,15,16], the high index differential-algebraic equations [17,18], the integro-differential equations [19] and the non-linear oscillators [20].

Recently, a new analytical technique is developed to solve fractional differential equations (FDEs) [21]. This technique, named as Fractional Differential Transform Method (FDTM), formulizes fractional power series likewise that DTM formulizes Taylor series. Studies were conducted about solutions of systems of fractional differential equations [22], systems of fractional partial differential equations [23], fractional-order integro-differential equations [24] and fractional differential-algebraic equations [25] using fractional differential transform method.

In this study, we give conformable fractional differential transform method (CFDTM) for conformable fractional derivative. CFDTM formulizes conformable fractional power series in a similar manner that FDTM formulizes fractional power series and DTM formulizes Taylor series.

## 2. Conformable Fractional Calculus

**Definition 2.1.** [1] Given a function $f:[0,\infty) \to \mathbb{R}$. Then the conformable fractional derivative of $f$ order $\alpha$ is defined by

$$(T_\alpha f)(t) = \lim_{\varepsilon \to 0} \frac{f(t + \varepsilon t^{1-\alpha}) - f(t)}{\varepsilon}$$

for all $t > 0, \alpha \in (0,1]$.

**Theorem 2.1.** [1] Let $\alpha \in (0,1]$ and $f, g$ be $\alpha$-differentiable at a point $t > 0$. Then
(1) $T_\alpha(af + bg) = a(T_\alpha f) + b(T_\alpha g)$, for all $a, b \in \mathbb{R}$
(2) $T_\alpha(t^p) = pt^{p-\alpha}$, for all $p \in \mathbb{R}$
(3) $T_\alpha(\lambda) = 0$, for all constant functions $f(t) = \lambda$
(4) $T_\alpha(fg) = f(T_\alpha g) + g(T_\alpha f)$
(5) $T_\alpha(f/g) = \frac{g(T_\alpha f) - f(T_\alpha g)}{g^2}$
(6) If, in addition, $f$ is differentiable, then $(T_\alpha f)(t) = t^{1-\alpha} \frac{df}{dt}(t)$.

**Definition 2.2.** [1] Given a function $f: [a, \infty) \to \mathbb{R}$. Then the conformable (left) fractional derivative of f order $\alpha$ is defined by
$$(T_\alpha^a f)(t) = \lim_{\varepsilon \to 0} \frac{f(t + \varepsilon(t-a)^{1-\alpha}) - f(t)}{\varepsilon}$$
for all $t > 0, \alpha \in (0,1]$.
All property in Theorem 2.1 is valid also for Definition 2.2 when $(t - a)$ is placed instead of $t$. Conformable fractional derivative of certain functions for Definition 2.2 is given following as:

(1) $T_\alpha^a((t-a)^p) = p(t-a)^{p-\alpha}$ for all $p \in \mathbb{R}$
(2) $T_\alpha^a \left( e^{\lambda \frac{(t-a)^\alpha}{\alpha}} \right) = \lambda e^{\lambda \frac{(t-a)^\alpha}{\alpha}}$
(3) $T_\alpha^a \left( \sin \left( \omega \frac{(t-a)^\alpha}{\alpha} + c \right) \right) = \omega \cos \left( \omega \frac{(t-a)^\alpha}{\alpha} + c \right), \ \omega, c \in \mathbb{R}$
(4) $T_\alpha^a \left( \cos \left( \omega \frac{(t-a)^\alpha}{\alpha} + c \right) \right) = -\omega \sin \left( \omega \frac{(t-a)^\alpha}{\alpha} + c \right), \ \omega, c \in \mathbb{R}$
(5) $T_\alpha^a \left( \frac{(t-a)^\alpha}{\alpha} \right) = 1$

**Definition 2.3.** Given a function $f: [a, \infty) \to \mathbb{R}$. Let $n < \alpha \leq n+1$ and $\beta = \alpha - n$. Then the conformable (left) fractional derivative of $f$ order $\alpha$, where $f^{(n)}(t)$ exists, is defined by
$$(T_\alpha^a f)(t) = (T_\beta^a f^{(n)})(t).$$

**Theorem 2.2.** [2] Assume $f$ is infinitely $\alpha$-differentiable function, for some $0 < \alpha \leq 1$ at a neighborhood of a point $t_0$. Then $f$ has the fractional power series expansion:
$$f(t) = \sum_{k=0}^\infty \frac{(T_\alpha^{t_0} f)^{(k)}(t_0)(t-t_0)^{\alpha k}}{\alpha^k k!}, \quad t_0 < t < t_0 + R^{1/\alpha}, \quad R > 0.$$
Here $(T_\alpha^{t_0} f)^{(k)}(t_0)$ denotes the application of the fractional derivative for $k$ times. For instance, let $y(t) = e^{\frac{t^\alpha}{\alpha}}$ and $t_0 = 0$, then
$$y(t) = e^{\frac{t^\alpha}{\alpha}} = \sum_{k=0}^\infty \frac{t^{k\alpha}}{\alpha^k k!}.$$
And for $y(t) = \frac{1}{1 - \frac{t^\alpha}{\alpha}}$, power series representation is

$$y(t) = \frac{1}{1-\frac{t^\alpha}{\alpha}} = \sum_{k=0}^{\infty} \frac{t^{\alpha k}}{\alpha^k}, \quad t \in [0,1).$$

## 3. Conformable Fractional Differential Transform Method

**Definition 3.1** Assume $f(t)$ is infinitely $\alpha$-differentiable function for some $\alpha \in (0,1]$. Conformable fractional differential transform of $f(t)$ is defined as

$$F_\alpha(k) = \frac{1}{\alpha^k k!} \left[ \left( T_\alpha^{t_0} f \right)^{(k)} (t) \right]_{t=t_0},$$

where $\left( T_\alpha^{t_0} f \right)^{(k)} (t)$ denotes the application of the fractional derivative $k$ times.

**Definition 3.2** Let $F_\alpha(k)$ be the conformable fractional differential transform of $f(t)$. Inverse conformable fractional differential transform of $F(k)$ is defined as

$$f(t) = \sum_{k=0}^{\infty} F_\alpha(k)(t-t_0)^{\alpha k} = \sum_{k=0}^{\infty} \frac{1}{\alpha^k k!} \left[ \left( T_\alpha^{t_0} f \right)^{(k)} (t) \right]_{t=t_0} (t-t_0)^{\alpha k}.$$

CFDT of initial conditions for integer order derivatives are defined as

$$F_\alpha(k) = \begin{cases} if\ \alpha k \in Z^+ & \frac{1}{(\alpha k)!} \left[ \frac{d^{\alpha k} f(t)}{dt^{\alpha k}} \right]_{t=t_0} \text{ for } k = 0,1,2,\ldots,\left(\frac{n}{\alpha}-1\right), \\ if\ \alpha k \notin Z^+ & 0 \end{cases}$$

where $n$ is the order of conformable fractional ordinary differential equation (CFODE).

**Theorem 3.1.** If $f(t) = u(t) \pm v(t)$, then $F_\alpha(k) = U_\alpha(k) \pm V_\alpha(k)$.

**Proof**

Conformable fractional differential transform of $u(t)$ and $v(t)$ can be written as the following:

$$U_\alpha(k) = \frac{1}{\alpha^k k!} \left[ \left( T_\alpha^{t_0} u \right)^{(k)} (t) \right]_{t=t_0}$$

$$V_\alpha(k) = \frac{1}{\alpha^k k!} \left[ \left( T_\alpha^{t_0} v \right)^{(k)} (t) \right]_{t=t_0}$$

Because of Theorem 2.2 (1), it is that

$$F_\alpha(k) = \frac{1}{\alpha^k k!} \left[ \left( T_\alpha^{t_0} (u(t) \pm v(t)) \right)^{(k)} (t) \right]_{t=t_0}$$

$$= \frac{1}{\alpha^k k!} \left[ \left( T_\alpha^{t_0} u \right)^{(k)} (t) \right]_{t=t_0} \pm \frac{1}{\alpha^k k!} \left[ \left( T_\alpha^{t_0} v \right)^{(k)} (t) \right]_{t=t_0} = U_\alpha(k) \pm V_\alpha(k)$$

□

**Theorem 3.2.** Let $c$ be a constant. If $f(t) = cu(t)$, then $F_\alpha(k) = cU_\alpha(k)$.

**Proof**

Conformable fractional differential transform of $u(t)$ is that

$$U_\alpha(k) = \frac{1}{\alpha^k k!} \left[ \left( T_\alpha^{t_0} u \right)^{(k)} (t) \right]_{t=t_0}$$

By the help of Theorem 2.2 (1), we can write

$$F_\alpha(k) = \frac{1}{\alpha^k k!} \left[ \left( T_\alpha^{t_0} cu \right)^{(k)} (t) \right]_{t=t_0} = \frac{c}{\alpha^k k!} \left[ \left( T_\alpha^{t_0} u \right)^{(k)} (t) \right]_{t=t_0} = cU_\alpha(k) \quad \square$$

**Theorem 3.3.** If $f(t) = u(t)v(t)$, then $F_\alpha(k) = \sum_{l=0}^{k} U_\alpha(l) V_\alpha(k-l)$.

**Proof**
By the help of Definition 3.1, $u(t)$ and $v(t)$ can be written that
$$u(t) = \sum_{k=0}^{\infty} U_\alpha(k)(t-t_0)^{\alpha k}, \quad v(t) = \sum_{k=0}^{\infty} V_\alpha(k)(t-t_0)^{\alpha k}.$$
Then, $f(t)$ is obtained as

$$f(x) = \sum_{k=0}^{\infty} U_\alpha(k)(t-t_0)^{\alpha k} \sum_{k=0}^{\infty} V_\alpha(k)(t-t_0)^{\alpha k}$$
$$= [U_\alpha(0) + U_\alpha(1)(t-t_0)^\alpha + U_\alpha(2)(t-t_0)^{2\alpha} + \cdots + U_\alpha(n)(t-t_0)^{n\alpha} + \cdots]$$
$$\times [V_\alpha(0) + V_\alpha(1)(t-t_0)^\alpha + V_\alpha(2)(t-t_0)^{2\alpha} + \cdots + V_\alpha(n)(t-t_0)^{n\alpha} + \cdots]$$
$$= [U_\alpha(0)V_\alpha(0) + (U_\alpha(0)V_\alpha(1) + U_\alpha(1)V_\alpha(0))(t-t_0)^\alpha$$
$$+ (U_\alpha(0)V_\alpha(2) + U_\alpha(1)V_\alpha(1) + U_\alpha(2)V_\alpha(0))(t-t_0)^{2\alpha} + \cdots]$$
$$= \sum_{k=0}^{\infty} \sum_{l=0}^{k} U_\alpha(l) V_\alpha(k-l) (t-t_0)^{k\alpha}$$

Hence, $F_\alpha(k)$ is found as
$$F_\alpha(k) = \sum_{l=0}^{k} U_\alpha(l) V_\alpha(k-l).$$

$\square$

On the other hand, if
$$f(t) = u_1(t) u_2(t) \ldots u_n(t),$$
it is clear
$$F_\alpha(k) = \sum_{k_{n-1}=0}^{k} \sum_{k_{n-2}}^{k_{n-1}} \cdots \sum_{k_2=0}^{k_3} \sum_{k_1=0}^{k_2} (U_1)_\alpha(k_1) (U_2)_\alpha(k_2 - k_1) \ldots (U_n)_\alpha(k - k_{n-1})$$

**Theorem 3.4.** If $f(t) = T_\alpha^{t_0}(u(t))$, then $F_\alpha(k) = \alpha(k+1) U_\alpha(k+1)$.

**Proof**
Let CFDT of $u(t)$ is as following
$$U_\alpha(k) = \frac{1}{\alpha^k k!} \left[ \left( T_\alpha^{t_0} u \right)^{(k)} (t) \right]_{t=t_0}$$
For $f(t) = T_\alpha^{t_0}(u(t))$,

$$F_\alpha(k) = \frac{1}{\alpha^k k!}\left[\left(T_\alpha^{t_0}(T_\alpha^{t_0}u)\right)^{(k)}(t)\right]_{t=t_0} = \frac{1}{\alpha^k k!}\left[\left(T_\alpha^{t_0}u\right)^{(k+1)}(t)\right]_{t=t_0}$$

$$= \alpha(k+1)\frac{1}{\alpha^{k+1}(k+1)!}\left[\left(T_\alpha^{t_0}u\right)^{(k+1)}(t)\right]_{t=t_0} = \alpha(k+1)U_\alpha(k+1)$$

□

**Theorem 3.5.** If $f(t) = T_\beta^{t_0}(u(t))$ for $m < \beta \leq m+1$, then
$F_\alpha(k) = U_\alpha(k+\beta/\alpha)\frac{\Gamma(k\alpha+\beta+1)}{\Gamma(k\alpha+\beta-m)}$.

**Proof**

Considering the initial conditions of the problem, we should seek the conformable differential transform of function $f(t)$. Therefore, we notice that

$$T_\beta^{t_0}(u(t)) = T_\beta^{t_0}\left(u(t) - \sum_{k=0}^{m}\frac{1}{k!}(t-t_0)^k u^{(k)}(0)\right).$$

Hence, if CFDT of $f(t) = T_\beta^{t_0}(u(t))$ is $F_\alpha(k)$, then CFDT of $T_\beta^{t_0}\left(u(t) - \sum_{k=0}^{m}\frac{1}{k!}(t-t_0)^k u^{(k)}(0)\right)$ is also $F_\alpha(k)$.

Let CFDT of $u(t)$ be $U_\alpha(k)$, then

$$f(t) = T_\beta^{t_0}\left(\sum_{k=0}^{\infty}U_\alpha(k)(t-t_0)^{k\alpha} - \sum_{k=0}^{m}\frac{1}{k!}(t-t_0)^k u^{(k)}(0)\right)$$

is written. Substituting $k\alpha$ in place of $k$ in the second series and considering the initial conditions, we obtain

$$f(t) = T_\beta^{t_0}\left(\sum_{k=0}^{\infty}U_\alpha(k)(t-t_0)^{k\alpha} - \sum_{k=0}^{\frac{\beta}{\alpha}-1}\frac{1}{(k\alpha)!}(t-t_0)^{k\alpha}u^{(k\alpha)}(0)\right)$$

$$= T_\beta^{t_0}\left(\sum_{k=0}^{\infty}U_\alpha(k)(t-t_0)^{k\alpha} - \sum_{k=0}^{\frac{\beta}{\alpha}-1}U_\alpha(k)(t-t_0)^{k\alpha}\right)$$

$$= T_\beta^{t_0}\left(\sum_{k=\frac{\beta}{\alpha}}^{\infty}U_\alpha(k)(t-t_0)^{k\alpha}\right).$$

$T_\beta^{t_0}\left(\sum_{k=\frac{\beta}{\alpha}}^{\infty}U_\alpha(k)(t-t_0)^{k\alpha}\right)$ is calculated as

$$T_\beta^{t_0}\left(\sum_{k=\frac{\beta}{\alpha}}^{\infty}U_\alpha(k)(t-t_0)^{k\alpha}\right) = \sum_{k=\frac{\beta}{\alpha}}^{\infty}U_\alpha(k)\frac{\Gamma(k\alpha+1)}{\Gamma(k\alpha-m)}(t-t_0)^{k\alpha}$$

$$= \sum_{k=0}^{\infty} U_\alpha(k + \beta/\alpha) \frac{\Gamma(k\alpha+\beta+1)}{\Gamma(k\alpha+\beta-m)} (t - t_0)^{k\alpha}.$$

According to this result, we can write
$$F_\alpha(k) = U_\alpha(k + \beta/\alpha) \frac{\Gamma(k\alpha+\beta+1)}{\Gamma(k\alpha+\beta-m)}. \quad \square$$

**Theorem 3.6.** If $f(x) = \left(T_{\beta_1}^{t_0} u_1\right)(t) \cdot \left(T_{\beta_2}^{t_0} u_2\right)(t) \ldots \left(T_{\beta_n}^{t_0} u_n\right)(t)$, then
$F_\alpha(k) =$
$\sum_{k_{n-1}=0}^{k} \sum_{k_{n-2}=0}^{k_{n-1}} \ldots \sum_{k_2=0}^{k_3} \sum_{k_1=0}^{k_2} \frac{\Gamma(k_1\alpha+\beta_1+1)}{\Gamma(k_1\alpha+\beta_1-m_1)} \frac{\Gamma((k_2-k_1)\alpha+\beta_2+1)}{\Gamma((k_2-k_1)\alpha+\beta_2-m_2)} \cdots \frac{\Gamma((k-k_{n-1})\alpha+\beta_n+1)}{\Gamma((k-k_{n-1})\alpha+\beta_n-m_n)} (U_1)_\alpha(k_1 + \beta_1/\alpha) (U_2)_\alpha(k_2 - k_1 + \beta_2/\alpha) \ldots (U_n)_\alpha(k - k_{n-1} + \beta_n/\alpha)$.
where $\frac{\beta_i}{\alpha} \in Z^+$ and $m_i < \beta_i \leq m_i + 1$ for $i = 1,2,\ldots,n$.

**Proof**

Let the conformable fractional differential transform of $\left(T_{\beta_i}^{t_0} u_i\right)(t)$ for $i = 1,2,\ldots,n$ be $(C_i)_\alpha(k)$. In this case, by the help of Theorem 3.3, we can write
$F_\alpha(k) = \sum_{k_{n-1}=0}^{k} \sum_{k_{n-2}=0}^{k_{n-1}} \ldots \sum_{k_2=0}^{k_3} \sum_{k_1=0}^{k_2} (C_1)_\alpha(k_1) (C_2)_\alpha(k_2 - k_1) \ldots (C_n)_\alpha(k - k_{n-1})$.
According to Theorem 3.4, it is that
$$(C_1)_\alpha(k_1) = \frac{\Gamma(k_1\alpha + \beta_1 + 1)}{\Gamma(k_1\alpha + \beta_1 - m_1)} (U_1)_\alpha(k_1 + \beta_1/\alpha)$$
$$(C_2)_\alpha(k_2 - k_1) = \frac{\Gamma((k_2 - k_1)\alpha + \beta_2 + 1)}{\Gamma((k_2 - k_1)\alpha + \beta_2 - m_2)} (U_2)_\alpha(k_2 - k_1 + \beta_2/\alpha)$$
$$\vdots$$
$(C_{n-1})_\alpha(k_{n-1} - k_{n-2})$
$$= \frac{\Gamma((k_{n-1} - k_{n-2})\alpha + \beta_{n-1} + 1)}{\Gamma((k_{n-1} - k_{n-2})\alpha + \beta_{n-1} - m_{n-1})} (U_{n-1})_\alpha(k_{n-1} - k_{n-2} + \beta_{n-1}/\alpha)$$
$$(C_n)_\alpha(k - k_{n-1}) = \frac{\Gamma((k-k_{n-1})\alpha+\beta_n+1)}{\Gamma((k-k_{n-1})\alpha+\beta_n-m_n)} (U_n)_\alpha(k - k_{n-1} + \beta_n/\alpha).$$
Hence, we obtain
$F_\alpha(k) =$
$\sum_{k_{n-1}=0}^{k} \sum_{k_{n-2}=0}^{k_{n-1}} \ldots \sum_{k_2=0}^{k_3} \sum_{k_1=0}^{k_2} \frac{\Gamma(k_1\alpha+\beta_1+1)}{\Gamma(k_1\alpha+\beta_1-m_1)} \frac{\Gamma((k_2-k_1)\alpha+\beta_2+1)}{\Gamma((k_2-k_1)\alpha+\beta_2-m_2)} \cdots \frac{\Gamma((k-k_{n-1})\alpha+\beta_n+1)}{\Gamma((k-k_{n-1})\alpha+\beta_n-m_n)} (U_1)_\alpha(k_1 + \beta_1/\alpha) (U_2)_\alpha(k_2 - k_1 + \beta_2/\alpha) \ldots (U_n)_\alpha(k - k_{n-1} + \beta_n/\alpha)$. $\square$

**Theorem 3.7.** If $f(t) = (t - t_0)^p$, then $F_\alpha(k) = \delta\left(k - \frac{p}{\alpha}\right)$ where
$$\delta(k) = \begin{cases} 1, & \text{if } k = 0 \\ 0, & \text{if } k \neq 0 \end{cases}.$$

**Proof**
CFDT of $f(t) = (t - t_0)^p$ is

$$F_\alpha(k) = \frac{1}{\alpha^k k!}\left[\left(T_\alpha^{t_0}((t-t_0)^p)\right)^{(k)}(t)\right]_{t=t_0}.$$

If the conformable fractional derivative of $f(t) = (t-t_0)^p$ is calculated for $k$ times, where $\alpha \in (0,1]$, then

$$= \frac{1}{\alpha^k k!}[p(p-\alpha)\ldots(p-(k-1)\alpha)(t-t_0)^{p-k\alpha}]_{t=t_0}$$

is obtained.

If $k = \frac{p}{\alpha}$, then

$$= \frac{1}{\alpha^{p/\alpha}\left(\frac{p}{\alpha}\right)!}\left(p(p-\alpha)\ldots(p-(p/\alpha - 1)\alpha)\right)$$

$$= \frac{1}{\alpha^{p/\alpha}\left(\frac{p}{\alpha}\right)!}\left(p(p-\alpha)\ldots(\alpha)\right)$$

$$= \frac{\alpha^{p/\alpha}}{\alpha^{p/\alpha}\left(\frac{p}{\alpha}\right)!}\left(\frac{p}{\alpha}\left(\frac{p}{\alpha}-1\right)\ldots(1)\right)$$

$$= \frac{\alpha^{p/\alpha}}{\alpha^{p/\alpha}\left(\frac{p}{\alpha}\right)!}\left(\frac{p}{\alpha}\right)!$$

$$= 1$$

Otherwise, for $t = t_0$, it is that

$$F_\alpha(k) = 0.$$

Hence

$$F_\alpha(k) = \delta\left(k - \frac{p}{\alpha}\right)$$

is obtained. where

$$\delta(k) = \begin{cases} 1, & if\ k = 0 \\ 0, & if\ k \neq 0 \end{cases}.$$

□

**Theorem 3.8.** If $f(t) = e^{\lambda\frac{(t-t_0)^\alpha}{\alpha}}$, then $F_\alpha(k) = \frac{\lambda^k}{\alpha^k k!}$, where $\lambda$ is constant.

**Proof**
CFDT of function $f$ is

$$F_\alpha(k) = \frac{1}{\alpha^k k!}\left[\left(T_\alpha^{t_0}\left(e^{\lambda\frac{(t-t_0)^\alpha}{\alpha}}\right)\right)^{(k)}(t)\right]_{t=t_0}.$$

Calculating the conformable fractional derivative of $f(t) = e^{\lambda\frac{(t-t_0)^\alpha}{\alpha}}$ for $k$ times, where $\alpha \in (0,1]$, we obtain

$$F_\alpha(k) = \frac{\lambda^k}{\alpha^k k!}.\ \square$$

**Theorem 3.9.** If $f(t) = \sin\left(\omega \frac{(t-t_0)^\alpha}{\alpha} + c\right)$, then $F_\alpha(k) = \frac{\omega^k}{\alpha^k k!} \sin\left(k\frac{\pi}{2} + c\right)$ and
If $f(t) = \cos\left(\omega \frac{(t-t_0)^\alpha}{\alpha} + c\right)$ for $\alpha \in (0,1)$, then $F_\alpha(k) = \frac{\omega^k}{\alpha^k k!} \cos\left(k\frac{\pi}{2} + c\right)$, where $w$ and $c$ are constant.

**Proof**

CFDT of $f(t) = \sin\left(\omega \frac{(t-t_0)^\alpha}{\alpha} + c\right)$ is written as

$$F_\alpha(k) = \frac{1}{\alpha^k k!}\left[\left(T_\alpha^{t_0}\left(\sin\left(\omega \frac{(t-t_0)^\alpha}{\alpha} + c\right)\right)\right)^{(k)}(t)\right]_{t=t_0}.$$

Thus, it is obtained that

$$F_\alpha(k) = \frac{\omega^k}{\alpha^k k!} \sin\left(k\frac{\pi}{2} + c\right).$$

Similarly, other part is proved. □

## 4. Applications

In this section, we give the solutions of some conformable fractional ordinary differential equations by the help of CFDTM. Additionally, we compare these solutions with the exact solutions.

**Example 1** Let's find the solution of the equation $y^{(\alpha)} + y = 0, y(0) = 1$ for $\alpha \in (0,1]$ by means of CFDTM. Exact solution of this equation is $y(t) = e^{-\frac{1}{\alpha}t^\alpha}$.
Because of Theorem 3.4, the equation above can be rewritten as following:
$$\alpha(k+1)Y_\alpha(k+1) + Y_\alpha(k) = 0, Y_\alpha(0) = 1.$$
Hence, recurrence relation is obtained as
$$Y_\alpha(k+1) = -\frac{1}{\alpha(k+1)} Y_\alpha(k), Y_\alpha(0) = 1.$$
For $k = 0,1,2,\ldots,n$,
$$Y_\alpha(1) = -\frac{1}{\alpha} Y_\alpha(0) = -\frac{1}{\alpha}$$
$$Y_\alpha(2) = -\frac{1}{2\alpha} Y_\alpha(1) = \frac{1}{2!\,\alpha^2}$$
$$Y_\alpha(3) = -\frac{1}{3\alpha} Y_\alpha(2) = -\frac{1}{3!\,\alpha^3}$$
$$\vdots$$
$$Y_\alpha(n) = \frac{(-1)^n}{n!\,\alpha^n}$$

is obtained. Hence the solution is

$$y(t) = \sum_{n=0}^\infty \frac{(-1)^n}{n!\,\alpha^n} t^{n\alpha} = e^{-\frac{1}{\alpha}t^\alpha}.$$

As seen, this solution is the same as the exact solution obtained previously.

**Example 2** We consider the fractional Riccati equation

$$y^{(\alpha)} = 1 + 2y + y^2$$

with the initial condition $y(0) = 0$. Exact solution of this equation is $y(t) = \frac{t^\alpha}{\alpha - t^\alpha}$.

Because of Theorem 3.2, Theorem 3.3, Theorem 3.4 and Theorem 3.6, we can write

$$\alpha(k+1)Y_\alpha(k+1) = \delta(k) - \sum_{l=0}^{k} Y_\alpha(l)Y_\alpha(k-l).$$

Thus, we obtain

$$Y_\alpha(k+1) = \frac{1}{\alpha(k+1)}\left(\delta(k) - \sum_{l=0}^{k} Y_\alpha(l)Y_\alpha(k-l)\right).$$

For $k = 0,1,2,\ldots$, the solution of this equation is

$$y(t) = \frac{t^\alpha}{\alpha} + \frac{t^{2\alpha}}{\alpha^2} + \frac{t^{3\alpha}}{\alpha^3} + \frac{t^{4\alpha}}{\alpha^4} + \frac{t^{5\alpha}}{\alpha^5} + \frac{t^{6\alpha}}{\alpha^6} + \frac{t^{7\alpha}}{\alpha^7} + \frac{t^{8\alpha}}{\alpha^8} + \frac{t^{9\alpha}}{\alpha^9} + \frac{t^{10\alpha}}{\alpha^{10}} + \cdots = \frac{t^\alpha}{\alpha}\left(1 + \frac{t^\alpha}{\alpha} + \left(\frac{t^\alpha}{\alpha}\right)^2 + \left(\frac{t^\alpha}{\alpha}\right)^3 + \left(\frac{t^\alpha}{\alpha}\right)^4 + \left(\frac{t^\alpha}{\alpha}\right)^5 + \left(\frac{t^\alpha}{\alpha}\right)^6 + \left(\frac{t^\alpha}{\alpha}\right)^7 + \left(\frac{t^\alpha}{\alpha}\right)^8 + \left(\frac{t^\alpha}{\alpha}\right)^9 + \cdots\right) = \frac{t^\alpha}{\alpha}\frac{1}{1-\frac{t^\alpha}{\alpha}} = \frac{t^\alpha}{\alpha - t^\alpha}.$$

This solution also is the same exact solution.

**Example 3** We consider the fractional Riccati equation

$$y^{(\alpha)} = 1 - y^2$$

with the initial condition $y(0) = 0$. $y(t) = \frac{e^{\frac{2}{\alpha}t^\alpha} - 1}{e^{\frac{2}{\alpha}t^\alpha} + 1}$ is the exact solution of this equation.

By the help of Theorem 3.3, Theorem 3.4 and Theorem 3.6, we can write

$$\alpha(k+1)Y_\alpha(k+1) = \delta(k) - \sum_{l=0}^{k} Y_\alpha(l)Y_\alpha(k-l).$$

Thereby, it is obtained that

$$Y_\alpha(k+1) = \frac{1}{\alpha(k+1)}\left(\delta(k) - \sum_{l=0}^{k} Y_\alpha(l)Y_\alpha(k-l)\right).$$

For $k = 0,1,2,\ldots$, the solution by means of CFDTM is found as

$$y(t) = \frac{t^\alpha}{\alpha} - \frac{t^{3\alpha}}{3\alpha^3} + \frac{2t^{5\alpha}}{15\alpha^5} - \frac{17t^{7\alpha}}{315\alpha^7} + \frac{62t^{9\alpha}}{2835\alpha^9} - \cdots.$$

The obtained solution of (1) above is the fractional power series expansion of the exact solution for the first ten terms.

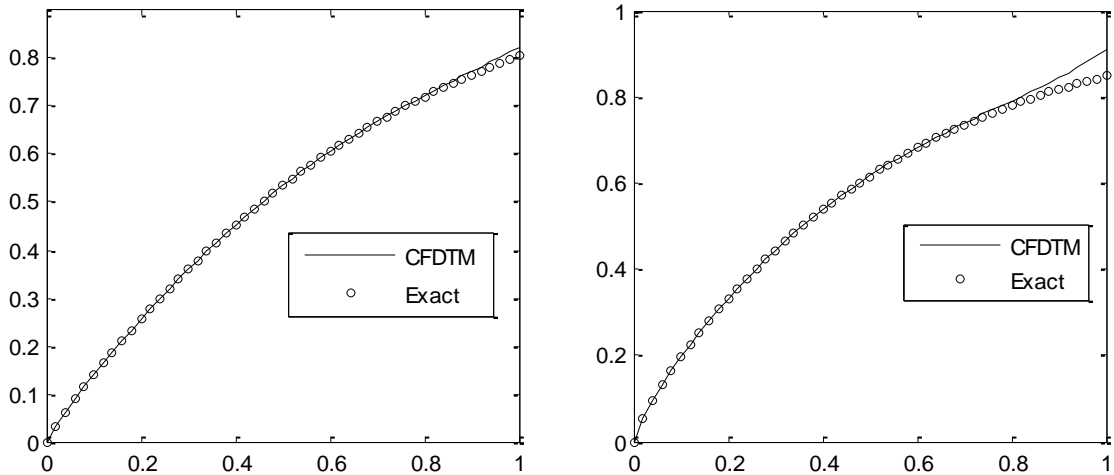

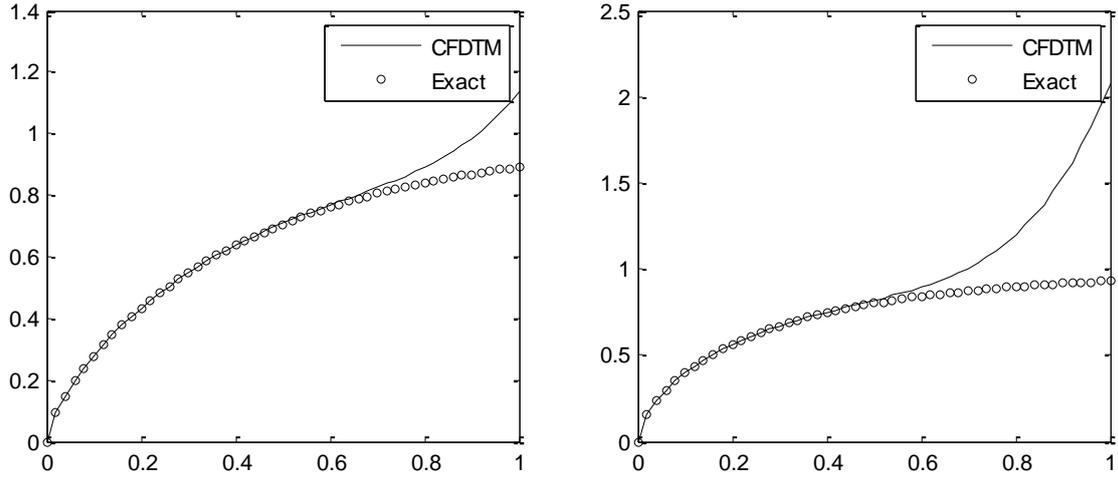

Fig. 1. The exact solutions versus the CFDTM solutions (solid lines) when $\alpha = 0.9, 0.8, 0.7, 0.6$ and $N = 10$ (the number of terms), respectively.

Figure 1 gives the relationship between the exact solutions and CFDTM solutions for Example 3. When the number of terms of series solution obtained with CFDTM is increased, the accuracy of the CFDTM solutions also increases. In this example, as the value of $\alpha$ decreases, CFDTM solutions retire from the exact solutions for interval [0,1], then a certain point.

**Example 4.** The equation Bagley-Torvik is following as:
$$A\frac{d^2y}{dt^2} + B\frac{d^{3/2}y}{dt^{3/2}} + Cy = f(t) \qquad (1)$$
with the boundary condition $y(0) = 1$ and $y'(0) = 1$. For $A = 1$, $B = 1$, $C = 1$ and $f(t) = 1 + t$, this equation is solved by the help of FDTM in [18]. In this work, selecting the order fraction as $\alpha = 0.5$, we solved the equation, with the same conditions, above via CFDTM. CFDT of boundary conditions are following as:
$$Y_\alpha(0) = 1, Y_\alpha(1) = 0, Y_\alpha(2) = 1 \text{ and } Y_\alpha(3) = 0.$$
Using Theorem 3.5 and Theorem 3.8, we can write CFDT of equation (1) as
$$Y_\alpha(k+4)\frac{\Gamma(0.5k+3)}{\Gamma(0.5k+1)} + Y_\alpha(k+3)\frac{\Gamma(0.5k+2.5)}{\Gamma(0.5k+0.5)} + Y_\alpha(k) = \delta(k) + \delta(k-2). \qquad (2)$$
Recurrence relation of (2) is
$$Y_\alpha(k+4) = \frac{\delta(k)+\delta(k-2)-Y_\alpha(k+3)(0.5k+0.5)(0.5k+1.5)-Y_\alpha(k)}{(0.5k+2)(0.5k+1)}.$$
Using the relation above, it is seen that
$$Y_\alpha(k+4) = 0, \ k = 0,1,2,...$$
Hence, the solution of (1) via CFDTM is found as $y(t) = 1 + t$ which is the same with the exact solution.

**Example 5.** Given a conformable fractional ordinary differential equation:
$$\frac{d^{3/2}y}{dt^{3/2}} - \frac{d^{1/2}y}{dt^{1/2}} = 0 \qquad (3)$$
with $y(0) = 0$ and $y'(0) = 1$. For $\alpha = 0.5$, by using theorem 3.6, CFDT of (3) is

$$Y_\alpha(k+3) = \frac{Y_\alpha(k+1)}{(0.5k+1.5)}.$$

For $k = 1, 2, \ldots$

$$Y_\alpha(4) = \frac{Y_\alpha(2)}{2} = \frac{1}{2}$$
$$Y_\alpha(5) = \frac{Y_\alpha(3)}{2.5} = 0$$
$$Y_\alpha(6) = \frac{Y_\alpha(4)}{3} = \frac{1}{2.3}$$
$$Y_\alpha(7) = \frac{Y_\alpha(5)}{3.5} = 0$$
$$\vdots$$

is obtained. Hence, the solution of (3) is found as

$$y(t) = t + \frac{1}{2!}t^2 + \frac{1}{3!}t^3 + \cdots$$
$$= \sum_{k=0}^{\infty} \frac{1}{k!} t^k - 1$$
$$= e^t - 1$$

which is the same of exact solution.

## 5. Conclusion

In this study, we present conformable fractional differential transform method (CFDTM) to find the numerical solution of conformable fractional ordinary differential equations. Then, we apply this new method to some conformable fractional ordinary differential equations. It is observed that CFDTM is an effective method for conformable fractional ordinary differential equations. Fractional power series solution is obtained with CFDTM. In some examples, the series solution obtained by the help of CFDTM can be written so as to exact solution. Otherwise, the number of terms in solution is increased to improve the accuracy of the obtained solution.